\documentclass[a4paper,11pt]{article}
 
\usepackage{latexsym}
\newtheorem{thm}{Theorem}[section]
\newtheorem{lem}[thm]{Lemma}

\newtheorem{prop}[thm]{Proposition}

\newtheorem{rem}[thm]{Remark}

\input epsf

\title{The Orchard relation of a generic symmetric or antisymmetric function}
\author{Roland Bacher}
\begin{document}
\maketitle

{\it Abstract: We associate to certain symmetric or antisymmetric
functions on the set ${E\choose d+1}$ of $(d+1)-$subsets in a finite
set $E$ an equivalence relation on $E$ and study some of its properties.}

\section{Definitions and main results}

We consider a finite set $E$ and denote by ${E\choose d}$ the set
of subsets containing exactly $d$ elements of $E$. In the sequel
we move often freely from sets to sequences: we identify a 
subset $\{x_1,\dots,x_d\}\in {E\choose d}$ with the finite
sequence $(x_1,\dots,x_d)$
where the order of the elements is for instance always
increasing with respect to a fixed total order on $E$.

A function $\varphi:{E\choose d}\longrightarrow {\mathbf R}$ 
is {\it symmetric} if
$$\varphi(x_1,\dots,x_i,x_{i+1},\dots,x_d)=\varphi(x_1\dots,x_{i-1},
x_{i+1},x_i,x_{i+1},\dots,x_d)$$
for $1\leq i<d$ and all $\{x_1,\dots,x_d\}\in {E\choose d}$.

Similarly, such a function 
$\varphi:{E\choose d}\longrightarrow{\mathbf R}$ is {\it antisymmetric} if
$$\varphi(x_1,\dots,x_i,x_{i+1},\dots,x_d)=-\varphi(x_1\dots,x_{i-1},
x_{i+1},x_i,x_{i+2},\dots,x_d)$$
for $1\leq i<d$ and all $x_1,\dots,x_d\in E$.

$\varphi$ is {\it generic} if $\varphi(x_1,\dots,x_d)\not= 0$
for all subsets $\{x_1,\dots,x_d\}\in {E\choose d}$
of $d$ distinct elements in $E$.

In the sequel of this paper all functions will be generic. We will
mainly be concerned with sign properties of generic 
symmetric or antisymmetric functions:
Given any symmetric generic function 
$\sigma:{E\choose d}\longrightarrow {\mathbf R}_{>0}$ and a symmetric or 
antisymmetric generic function $\varphi:{E\choose d}\longrightarrow
{\mathbf R}^d$, the two functions 
$$(x_1,\dots,x_d)\longmapsto \varphi(x_1,\dots,x_d)$$
and
$$(x_1,\dots,x_d)\longmapsto \sigma(x_1,\dots,x_d)\varphi(x_1,\dots,x_d)$$
behave similarly with respect to all properties adressed in this paper.

We have also an obvious sign rule: symmetric or antisymmetric
functions on ${E\choose d}$ behave with respect to multiplication
like the elements of the multiplicative group $\{\pm 1\}$ 
with symmetric functions corresponding to $1$ and antisymmetric 
functions corresponding to $-1$.

We fix now a generic symmetric or antisymmetric function 
$\varphi:{E\choose d+1}
\longrightarrow{\mathbf R}$. Consider two elements
$a,b\in E$. A subset
$\{x_1,\dots,x_d\}\in {E\setminus\{a,b\}\choose d}$ not containing $a$ and $b$ 
{\it separates} $a$ from $b$ with respect to $\varphi$ if 
$$\varphi(x_1,\dots,x_d,a)\ \varphi(x_1,\dots,x_d,b)<0$$
(this definition is of course independent of the particular linear
order $(x_1,\dots,x_d)$ on the set $\{x_1,\dots,x_d\}$).

We denote by $n(a,b)=n_\varphi(a,b)$ the number of subsets in 
${E\setminus\{a,b\}\choose d}$
separating $a$ from $b$ (with respect to the function $\varphi$).

\begin{prop} \label{propabc}
(i) If $\varphi:{E\choose d+1}\longrightarrow {\mathbf R}$ is
symmetric and generic then
$$n(a,b)+n(b,c)+n(a,c)\equiv 0\pmod 2$$
for any subset $\{a,b,c\}$ of $3$ distinct elements in $E$.

\ \ (ii) If $\varphi:{E\choose d+1}\longrightarrow {\mathbf R}$ is
antisymmetric and generic then
$$n(a,b)+n(b,c)+n(a,c)\equiv {\sharp(E)-3\choose d-1}\pmod 2\ $$
for any subset $\{a,b,c\}$ of $3$ distinct elements in $E$.
\end{prop}

{\bf Proof.} Consider first a subset $\{x_1,\dots,x_d\}$ not intersecting
$\{a,b,c\}$. Such a subset separates no pair of elements
in $\{a,b,c\}$ if 
$$\varphi(x_1,\dots,x_d,a),\ \varphi(x_1,\dots,x_d,b)\hbox{ and }
\varphi(x_1,\dots,x_d,c)$$
all have the same sign. Otherwise, consider a reordering 
$\{a',b',c'\}=\{a,b,c\}$ such that 
$\varphi(x_1,\dots,x_d,a')\ \varphi(x_1,\dots,x_d,b')<0$ and 
$\varphi(x_1,\dots,x_d,a')\ \varphi(x_1,\dots,x_d,c')<0$. The subset 
$\{x_1,\dots,x_d\}$ contributes in this case $1$ to $n(a',b'),\ n(a',c')$
and $0$ to $n(b',c')$. Such a subset
$\{x_1,\dots,x_d\}\in {E\setminus\{a,b,c\}\choose d}$ 
yields hence always an even contribution (0 or 2)
to the sum $n(a,b)+n(a,c)+n(b,c)$.

Consider now a subset $\{x_1,\dots,x_{d-1}\}\in 
{E\setminus\{a,b,c\}\choose d-1}$.
We have to understand the contributions of the sets
$$\begin{array}{lll}
\{x_1,\dots,x_{d-1},c\}\quad &\hbox{to }\ &n(a,b)\ ,\\
\{x_1,\dots,x_{d-1},b\}\quad &\hbox{to }\ &n(a,c)\ ,\\
\{x_1,\dots,x_{d-1},a\}\quad &\hbox{to }\ &n(b,c)\ .
\end{array}$$
Since the product of the six factors
$$\begin{array}{c}
\varphi(x_1,\dots,x_{d-1},c,a)\ \varphi(x_1,\dots,x_{d-1},c,b)\\
\varphi(x_1,\dots,x_{d-1},b,a)\ \varphi(x_1,\dots,x_{d-1},b,c)\\
\varphi(x_1,\dots,x_{d-1},a,b)\ \varphi(x_1,\dots,x_{d-1},a,c)
\end{array}$$
is always positive (respectively negative) for a generic symmetric
(respectively antisymmetric) function, such a subset 
$\{x_1,\dots,x_{d-1}\}$ yields an even contribution to 
$n(a,b)+n(a,c)+n(b,c)$ in the symmetric case and an odd contribution
in the antisymmetric case. 

Proposition \ref{propabc} follows now from the fact that ${E\setminus\{a,b,c\}\choose d-1}$ has ${\sharp(E)-3\choose d-1}$ elements.\hfill $\Box$

Given a generic symmetric or antisymmetric function
$\varphi:{E\choose d+1}\longrightarrow {\mathbf R}$ on some finite set $E$
we set $x\sim y$ if either $x=y\in E$ or if
$$n(x,y)\equiv 0\pmod 2\qquad \hbox{ for symmetric $\varphi$}$$
respectively
$$n(x,y)\equiv {\sharp(E)-3\choose d-1}\pmod 2\qquad \hbox{ for antisymmetric 
$\varphi$}\ .$$ We call the relation $\sim$ defined in this way
on the set $E$ the {\it Orchard relation}.

\begin{thm} The Orchard relation is an equivalence relation having at most two
classes.
\end{thm}

{\bf Proof.} Reflexivity and symmetry are obvious. Transitivity 
follows easily from Proposition \ref{propabc}.

If $a\not\sim b$ and $b\not\sim c$ then $n(a,b)+n(b,c)$ is even. It follows
then from Proposition \ref{propabc} that $a\sim c$. \hfill $\Box$

{\bf Example.} A {\it tournament} is a generic antisymmetric function
${E\choose 1+1}\longrightarrow \{\pm 1\}$. It encodes for instance
orientations of all edges in the complete graph with vertices $E$ 
and can be summarized by an
antisymmetric matrix $A$ with coefficients in $\{\pm 1\}$.

Given such a matrix $A$ with coefficients $a_{i,j},\ 1\leq i,j\leq n$,
we have
$$n_A(i,j)=\frac{n-2-\sum_k a_{ik}a_{jk}}{2}\ .$$
This implies $i\sim_A j$ if and only if
$$\sum_ka_{ik}a_{jk}\equiv n\pmod 4$$
for $i\not= j$. In the language of tournaments (cf. for instance \cite{M}),
this result can be restated in terms of score vectors: Two elements $i$ and
$j$ are Orchard equivalent if and only if the corresponding
coefficients of the score vector (counting the number of $1$'s in line $i$
respectively $j$) have the same parities.

{\bf Main Example.} A finite
set ${\cal P}=\{P_1,\dots,P_n\}\subset {\mathbf R}^d$ of 
$n>d$ points in real affine space
${\mathbf R}^d$
is {\it generic} if the affine span of any  subset containing $(d+1)$ points 
in $\cal P$ is all of ${\mathbf R}^{d}$. Such a generic set $\cal P$
is endowed with a generic antisymmetric function by restricting
$$\varphi(x_0,\dots,x_d)=\det(x_1-x_0,x_2-x_0,\dots,x_d-x_0)$$
to ${{\mathcal P}\choose d+1}$. The Orchard relation partitions hence a 
generic subset ${\cal P}\subset {\mathbf R}^d$ into two (generally
non-empty) subsets. Its name originates from the fact that
the planar case ($d=2$) yields a 
natural rule to plant trees of two different species at specified
generic locations in an orchard, see \cite{B1} and \cite{BG}.

\begin{prop} \label{flip} Given a finite set $E$ 
let $\varphi$ and $\psi$ be two generic symmetric or antisymmetric
functions on ${E\choose d+1}$.

\ \ (i) If the numbers 
$$\varphi(x_0,\dots,x_d)\ \psi(x_0,\dots,x_d)$$
have the same sign for all $\{x_0,\dots,x_d\}\in{E\choose d+1}$ then the two Orchard relations
$\sim_\varphi$ and $\sim_\psi$ induced by $\varphi$ and $\psi$ 
coincide.

\ \ (ii) If there exists exactly one subset ${\cal F}=\{x_0,\dots,x_d\}\in 
{E\choose d+1}$ such that  
$$\varphi(x_0,\dots,x_d)\ \psi(x_0,\dots,x_d)<0$$
then the restrictions of $\sim_\varphi$ and $\sim_\psi$ to the
two subsets
${\cal F}$ and $E\setminus {\cal F}$ coincide but
$a\sim_\varphi b\Longleftrightarrow a\not\sim_\psi b$ for $a\in{\cal F}$ 
and $b\in E\setminus{\cal F}$.
\end{prop} 

We call two symmetric or antisymmetric functions $\varphi$ and $\psi$
satisfying the condition of assertion (ii) above {\it flip-related}.
Coulouring the equivalence classes of an Orchard relation with two 
distinct coulours, one can express assertion (ii) by the statement that 
changing a generic (symmetric or antisymmetric) function by a flip 
switches the coulours in the flip-set ${\cal F}=\{x_0,\dots,x_d\}$
and leaves the coulours of the remaining elements unchanged. 

Assertion (i) shows that we can restrict our attention to symmetric or
antisymmetric functions from ${E\choose d+1}$ into
$\{\pm 1\}$ when studying properties of the Orchard relation.

{\bf Proof of Proposition \ref{flip}.} Assertion (i) is obvious.

For proving assertion (ii) it is enough to remark that the numbers 
$n_\varphi(a,b)$ and $n_\psi(a,b)$ of separating sets (with
respect to $\varphi$ and $\psi$) are identical if either $\{a,b\}\subset
{\cal F}$ or $\{a,b\}\subset E\setminus{\cal F}$ and
they differ by exactly one in the remaining cases. \hfill $\Box$

\section{An easy characterisation in the symmetric case}

In this section we give a different and rather 
trivial description of the Orchard relation
in the symmetric case.

Given a generic symmetric function $\varphi:{E\choose d+1}
\longrightarrow {\mathbf R}$ on some finite set $E$ we consider the function
$$\mu(x)=\sharp(\{\{x_1,\dots,x_d\}\in{E\setminus \{x\}\choose d}\ \vert
\ \varphi(x,x_1,\dots,x_d)>0\})$$
from $E$ to ${\mathbf N}$.

\begin{thm} Two elements $x,y\in E$ are Orchard equivalent with respect to 
$\varphi$ if and only if $\mu(x)\equiv \mu(y)\pmod 2$.
\end{thm}

{\bf Proof.} The result holds if $\varphi$ is the constant function
$$\varphi(x_0,\dots,x_d)=1$$
for all $\{x_0,\dots,x_d\}\in {E\choose d+1}$. 

Given two generic symmetric functions $\varphi,\psi$
related by a flip with respect to the 
set ${\cal F}=\{x_0,\dots,x_d\}\in {E\choose
d+1}$ we have
$$\mu_\varphi(x)=\mu_\psi(x)$$ 
if $x\not\in {\cal F}$ and 
$$\mu_\varphi(x)=\mu_\psi(x)\pm 1$$ 
otherwise. Proposition \ref{flip} 
implies hence the result since any generic symmetric
function can be related by a finite number of flips to the constant
function.\hfill $\Box$

\section{Reducing $d$}

Let $\varphi:{E\choose d+1}\longrightarrow {\mathbf R}$ be a generic
symmetric or antisymmetric function. Consider the function
$$R\varphi:{E\choose d}\longrightarrow {\mathbf R}$$
defined by
$$R\varphi(x_1,\dots,x_d)=\prod_{x\in E\setminus\{x_1,\dots,x_d\}}
\varphi(x,x_1,\dots,x_d)\ .$$

$R\varphi$ is generic symmetric if $\varphi$ is generic symmetric.

For $\varphi$ generic antisymmetric, the function $R\varphi$ is
generic symmetric if $\sharp(E)\equiv d\pmod 2$ and $R\varphi$
is generic antisymmetric otherwise.

Dependencies of the Orchard relations
associated to $\varphi$ and $R\varphi$ are described by the following result.

\begin{prop} \label{proprsym}
Let $\varphi:{E\choose d+1}\longrightarrow {\mathbf R}$ be a 
generic symmetric or antisymmetric function.

\ \ (i) If $d\equiv 0\pmod 2$ then the Orchard relation of $R\varphi$ is 
trivial (i.e. $x\sim_{R\varphi} y$ for all $x,y\in E$).

\ \ (ii) If $d\equiv 1\pmod 2$ then the Orchard relations $\sim_\varphi$ and
$\sim_{R\varphi}$ coincide on $E$.
\end{prop}

The main ingredient of the proof is the following lemma.

\begin{lem} \label{Rflip}
Let $\varphi,\psi:{E\choose d+1}\longrightarrow {\mathbf R}$ 
be two generic symmetric or 
antisymmetric functions which are flip-related with respect to
the set ${\cal F}=\{x_0,\dots,x_d\}\in {E\choose d+1}$. Then 
$$R\varphi(y_1,\dots,y_d)\  R\psi(y_1,\dots,y_d)<0$$
if $\{y_1,\dots,y_d\}\subset {\cal F}$ and 
$$R\varphi(y_1,\dots,y_d)\ R\psi(y_1,\dots,y_d)>0$$
otherwise.
\end{lem}

{\bf Proof of Lemma \ref{Rflip}.} If $\{y_1,\dots,y_d\}\not\subset 
{\cal F}$ then
$\varphi(x,y_1,\dots,y_d)=\psi(x,y_1,\dots,y_d)$ for all 
$x\in E\setminus\{y_1,\dots,y_d\}$ and hence $R\varphi(y_1,\dots,y_d)=
R\psi(y_1,\dots,y_d)$. Otherwise, exactly one factor of the product
yielding $R\psi(y_1,\dots,y_d)$ changes sign with respect to
the factors yielding $R\varphi(y_1,\dots,y_d)$.\hfill $\Box$

{\bf Proof of Proposition \ref{proprsym}.} We consider first the case where
$\varphi:{E\choose d+1}\longrightarrow {\mathbf R}$ is generic and symmetric.

Proposition \ref{proprsym}  holds then for the constant symmetric application
$\varphi:{E\choose d+1}\longrightarrow \{1\}$.

Two generic symmetric functions $\varphi,\psi$ on ${E\choose d+1}$ which
are flip-related with respect to ${\cal F}=\{x_0,\dots,x_d\}$ give rise to 
$R\varphi$ and $R\psi$ which are related through $d+1$ flips with respect to all $d+1$ elements in ${{\cal F}\choose d}$ by Lemma \ref{Rflip}. 
Proposition \ref{flip} implies hence the result since an element of $E\setminus {\cal F}$ is contained in no element of ${{\cal F}\choose d}$ and since 
all elements of ${\cal F}$ are contained in exactly $d$ such sets.

Second case: $\varphi:{E\choose d+1}\longrightarrow {\mathbf R}$ generic and
antisymmetric. This case is slightly more involved. As in the symmetric
case, we prove the result for a particular function $\varphi$ and use
the fact that flips of $\varphi$ affect the Orchard relation $\sim_{R\varphi}$
only for odd $d$. This shows that it is enough to prove that 
$\sim_{R\varphi}$ is trivial for a particular function $\varphi$ in the case
of even $d$ and that $\sim_{R\varphi}$ and $\sim_\varphi$ coincide (for
a particular generic antisymmetric function $\varphi$) in the case of 
odd $d$.

We consider now the set $E=\{1,\dots,n\}$ endowed with the generic
antisymmetric function $\varphi:{E\choose d+1}\longrightarrow\{\pm 1\}$
defined by
$$\varphi(i_0,\dots,i_d)=1$$
for all $1\leq i_0<i_1<\dots<i_d\leq n$.

Each element of ${E\setminus \{i,i+1\}\choose d-1}$ separates then
$i$ from $i+1$ with respect to the generic function $R\varphi$. We have
indeed
$$\begin{array}{l}
\displaystyle R\varphi(j_1,\dots,j_{d-1},i)\\
\displaystyle \quad =\varphi(i+1,j_1,\dots,j_{d-1},i)\prod_{j\in E\setminus \{j_1,\dots,j_{d-1},i,
i+1\}}\varphi(j,j_1,\dots,j_{d-1},i)\\
\displaystyle \quad
=-\varphi(i,j_1,\dots,j_{d-1},i+1)\prod_{j\in E\setminus \{j_1,\dots,j_{d-1},i,
i+1\}}\varphi(j,j_1,\dots,j_{d-1},i+1)\\
\displaystyle \qquad =-R\varphi(j_1,\dots,j_{d-1},i+1)\end{array}$$
showing that the number $n_{R\varphi}(i,i+1)$ of sets separating $i$ from 
$i+1$ equals ${n-2\choose d-1}$.

The proof splits now into four cases according to the parities of 
$n$ and $d$.

If $n\equiv d\equiv 0\pmod 2$, then $R\varphi$ is symmetric and 
${n-2\choose d-1}$ is even (recall that 
$${\sum_{i=0} \nu_i2^i\choose \sum_{i=0}\kappa_i 2^i}\equiv \prod_i {\nu_i\choose \kappa_i}\pmod 2$$
for $\nu_i,\kappa_i\in \{0,1\}$, cf. for instance 
Exercice 5.36 in Chapter 5 of \cite{GKP}). 
Since $n_{R\varphi}(i,i+1)={n-2\choose d-1}$ is even for all $i<n$,
the Orchard relation $\sim_{R\varphi}$ associated to the 
symmetric function $R\varphi$ is trivial.

If $n\equiv 1\pmod 2,\ d\equiv 0\pmod 2$, then $R\varphi$ is antisymmetric.
We have then ${n-3\choose d-1}\equiv 0\pmod 2$
and thus ${n-3\choose d-2}\equiv {n-3\choose d-2}+{n-3\choose d-1}=
{n-2\choose d-1}\pmod 2$ which implies
again the triviality of the Orchard relation $\sim_{R\varphi}$ since we have
$n_{R\varphi}(i,i+1)={n-2\choose d-1}\equiv {n-3\choose d-2}\pmod 2$
which shows $i\sim_{R\varphi} (i+1)$ for all $i$.

If $n\equiv d\equiv 1\pmod 2$ then $R\varphi$ is symmetric. Since 
${n-3\choose d-2}\equiv 0\pmod 2$ we have
${n-2\choose d-1}={n-3\choose d-2}+{n-3\choose d-1}
\equiv {n-3\choose d-1}\pmod 2$ proving
that the Orchard relations $\sim_\varphi$ and $\sim_{R\varphi}$
coincide.

If $n\equiv 0\pmod 2,\ d\equiv 1\pmod 2$, then $R\varphi$ is antisymmetric.
The equality ${n-2\choose d-1}={n-3\choose d-2}+{n-3\choose d-1}$
implies ${n-3\choose d-1}\equiv {n-3\choose d-2}+{n-2\choose d-1}\pmod 2$.
This shows hat the Orchard relations $\sim\varphi$ and $\sim_{R\varphi}$
coincide.\hfill $\Box$

\section{Homology}

We recall that $R\varphi:{E\choose d}\longrightarrow {\mathbf R}$ 
is defined by
$$R\varphi(x_1,\dots,x_d)=\prod_{x\in E\setminus\{x_1,\dots,x_d\}}
\varphi(x,x_1,\dots,x_d)$$
for a given generic symmetric or antisymmetric function
$\varphi:{E\choose d+1}\longrightarrow {\mathbf R}$.

\begin{lem} We have 
$$R(R\varphi)(x_1,\dots,x_{d-1})\in \epsilon^{\sharp(E)-d+1\choose 2}
{\mathbf R}_{>0}$$
where $\epsilon=1$ if $\varphi$ is generic and symmetric and $\epsilon=-1$ if
$\varphi$ is generic and antisymmetric.
\end{lem}

{\bf Proof.} Setting ${\cal S}=\{x_1,\dots,x_{d-1}\}$ we have
$$\begin{array}{l}
R(R\varphi)(x_1,\dots,x_{d-1})\\
\quad =\prod_{y\in E\setminus{\cal S}}R\varphi(y,x_1,\dots,x_{d-1})=
\prod_{x\not=y\in E\setminus{\cal S}}\varphi(x,y,x_1,\dots,x_{d-1})\\
\quad =\prod_{\{x,y\}\in {E\setminus {\cal S}\choose 2}}
\varphi(x,y,x_1,\dots,x_{d-1})\varphi(y,x,x_1,\dots,x_{d-1})
\end{array}$$
which is positive if $\varphi$ is symmetric or if ${\sharp(E)-d+1\choose 2}$
is even and negative otherwise.\hfill $\Box$

Writing as in the beginning $[n]=\{1,\dots,n\}$,
the set $\{\pm 1\}^{[n]\choose d+1}$ (endowed with the the usual
product of functions) of all symmetric generic functions
${[n]\choose d+1}\longrightarrow \{\pm 1\}$ is a vector
space of dimension ${n\choose d+1}$ over the field ${\mathbf F}_2$ of
2 elements. The map $R$ considered above defines group homomorphisms
between these vector spaces and the above Lemma allows to define 
homology groups. These groups are however all trivial except for
$d=0$ since one obtains the ordinary 
(simplicial) homology with coefficients in ${\mathbf F}_2$ 
of an $(n-1)$ dimensional simplex.

\section{Increasing $d$}

This section is a close analogue of section 3.

Given a generic symmetric or antisymmetric function 
${E\choose d+1}\longrightarrow {\mathbf R}$ we define a function
$A\varphi:{E\choose d+2}\longrightarrow {\mathbf R}$ by setting
$$A\varphi(x_0,\dots,x_{d+1})=\prod_{i=0}^{d+1}\varphi(x_0,\dots,x_{i-1},
x_{i+1},\dots,x_{d+1})\ .$$

The function $A\varphi$ is generic symmetric if $\varphi$ is symmetric.
For $\varphi$ antisymmetric it is generic symmetric if $d\equiv 0\pmod 2$ 
and generic antisymmetric otherwise. 

The dependency between the Orchard relations 
$\sim_\varphi$ and $\sim_{A\varphi}$ for
a generic symmetric or antisymmetric function $\varphi:{E\choose d+1}
\longrightarrow {\mathbf R}$ is described
by the following result.

\begin{prop} \label{Aprop} Let $\varphi:{E\choose d+1}
\longrightarrow {\mathbf R}$ be a generic symmetric or antisymmetric function.

The Orchard relation $\sim_{A\varphi}$ of $A\varphi$ is trivial if 
$\sharp(E)\equiv d\pmod 2$. Otherwise, the Orchard relations 
$\sim_{\varphi}$ and $\sim_{A\varphi}$ of $\varphi$ and $A\varphi$
coincide.
\end{prop}

The main ingredient of the proof is the following lemma whose easy proof
is left to the reader.

\begin{lem} \label{Alem} Let 
$\varphi,\psi:{E\choose d+1}\longrightarrow {\bf R}$ be two generic symmetric
or antisymmetric functions which are flip-related with respect to the set
${\cal F}=\{x_0,\dots,x_d\}$. Then

$$A\varphi(y_0,\dots,y_{d+1})\ A\psi(y_0,\dots,y_{d+1})>0$$
if ${\cal F}\not\subset \{y_0,\dots,y_{d+1}\}$ and 
$$A\varphi(y_0,\dots,y_{d+1})\ A\psi(y_0,\dots,y_{d+1})<0$$
otherwise.
\end{lem}

{\bf Proof of Proposition \ref{Aprop}.}
Lemma \ref{Alem} shows that $\sim_{A\varphi}$ is independent of
$\varphi$ if $\sharp(E)\equiv d\pmod 2$. Otherwise, the Orchard relations of 
$\varphi$ and $A\varphi$ behave in a similar way under flips. 
Indeed, given $\psi$ which is flip-related with 
flipset ${\cal F}=\{x_0,\dots,x_d\}$ to $\varphi$ the functions 
$A\psi$ and $A\varphi$ are related through $(\sharp(E)-(d+1))$ flips with 
flipsets ${\cal F}\cup\{x\}$, $x\in E\setminus {\cal F}$. 
Each element of $E\setminus {\cal F}$ is hence flipped once and each 
element of ${\cal F}$ is flipped $\sharp (E)-(d+1)$ times. 

Proposition \ref{flip} implies hence that $\sim_{A\varphi}$ is
independent of $\varphi$ if $1\equiv \sharp(E)-(d+1)$ and that 
$\sim_\varphi$ and $\sim_{A\varphi}$ behave similarly under flips 
otherwise. It is hence enough to proof Proposition \ref{Aprop}
in a particular case. 

If $\varphi$ is symmetric, then Proposition \ref{Aprop} clearly holds
for the constant application $\varphi:{E\choose d+1}\longrightarrow \{1\}$.

In the antisymmetric case we set $E=\{1,\dots,n\}$
and we consider the generic
antisymmetric function $\varphi:{E\choose d+1}\longrightarrow\{\pm 1\}$
defined by
$$\varphi(i_0,\dots,i_d)=1$$
for all $1\leq i_0<i_1<\dots<i_d\leq n$. The function $A\varphi:{E\choose
d+2}\longrightarrow \{\pm 1\}$ is now given by
$$A\varphi(i_0,\dots,i_d,i_{d+1})=1$$
for all $1\leq i_0<i_1<\dots<i_d<i_{d+1}\leq n$.
The numbers $n_{A\varphi}(i,i+1)$ of subsets in ${E\setminus \{i,i+1\}\choose d+1}$ separating $i$ from $i+1$ are hence all $0$ and we split the discussion
into several cases according to the parities of $n=\sharp(E)$ and $d$.

$n\equiv d\equiv 0\pmod 2$ implies $A\varphi$ symmetric and hence 
$\sim_{A\varphi}$ trivial. 

$n\equiv 1\pmod 2,\ d\equiv 0\pmod 2$ implies $A\varphi$ symmetric
and hence $\sim_{A\varphi}$ trivial. Since then ${n-3\choose d-1}\equiv 
0\pmod 2$ we have also $\sim_\varphi$ trivial.

$n\equiv 0\pmod 2,\ d\equiv 1\pmod 2$ implies $A\varphi$ antisymmetric.
We have then ${n-3\choose d-1}+{n-3\choose d}={n-2\choose d}\equiv 0\pmod 2$
proving equality of the two Orchard relations $\sim_{\varphi}$ and
$\sim_{A\varphi}$.

$n\equiv d\equiv 1\pmod 2$ implies $A\varphi$ antisymmetric and ${n-3\choose d}
\equiv 0\pmod 2$ thus proving triviality of the Orchard relation 
$\sim_{A\varphi}$. \hfill $\Box$

\begin{rem} One sees easily that the function $A(A\varphi)$
is strictly positive for a generic symmetric or antisymmetric
function $\varphi:{E\choose d+1}\longrightarrow {\mathbf R}$.

This allows the definition of cohomology groups on the set
of generic symmetric functions ${E\choose d+1}\longrightarrow \{\pm 1\}$.
The resulting groups are of course not interesting since
this boils down once more to the cohomology groups of the
$(\sharp(E)-1)-$dimensional simplex with 
coefficients in the field ${\mathbf F}_2$ of $2$ elements.
\end{rem}


\begin{thebibliography}{99}


\bibitem{B1} R. Bacher, 
{\it An Orchard Theorem}, preprint CO/0206266, 14 pages.

\bibitem{BG} R. Bacher, D. Garber, {\it Chromatic Properties of generic
planar Configurations of Points}, preprint GT/0210051, 19 pages.

\bibitem{GKP} Graham, Knuth, 0. Patashnik, {\it Concrete Mathematics: A foundation for Computer Science}, 2nd ed, Addison-Wesely (1994).

\bibitem{M} J.W. Moon, {\it Topics on Tournaments}, Holt, Rinehart and Winston
(1968).

\end{thebibliography}
\end{document}